\documentstyle[12pt]{article}%
\newcommand{\mt}[1]{\mathop{\rm #1}\nolimits}%
\newcommand{\mod}{\mathop{\rm mod}\nolimits}%
\newcommand{\wt}{\mathop{\rm wt}\nolimits}%
\newcommand{\pal}{\text{---}}%

\newcommand{\bb}[1]{{\BBB #1}}%
\newcommand{\CC}{{\BBB C}}%
\newcommand{\ZZ}{{\BBB Z}}%
\newcommand{\QQ}{{\BBB Q}}%
\newcommand{\PP}{{\BBB P}}%
\newcommand{\NN}{{\BBB N}}%
\newcommand{\OOO}{{\cal O}}%

\newcommand{\text}[1]{\mbox{{\rm #1}}}%
\newcommand{\qq}{$\ \Box$}%
\newcommand{\cyc}[1]{\bb{Z}_{#1}}%
\newcommand{\3}{^{\sharp}}

\font\twelvegtc=eufm10 scaled 1200%
\font\ninegtc=eufm9%
\font\sevengtc=eufm7%
\newfam\gtcfam%
\textfont\gtcfam=\twelvegtc%
\scriptfont\gtcfam=\ninegtc%
\scriptscriptfont\gtcfam=\sevengtc%
\font\twelveBBB=msbm10 scaled 1200%
\font\tenBBB=msbm10%
\font\sevenBBB=msbm7%
\newfam\BBBfam%
\def\BBB{\fam\BBBfam\twelveBBB}%
\textfont\BBBfam=\twelveBBB%
\scriptfont\BBBfam=\tenBBB%
\scriptscriptfont\BBBfam=\sevenBBB%

\newcounter{THNO}[section]%
\renewcommand{\theTHNO}{(\arabic{section}.\arabic{THNO}) }%

\newcounter{sub}[THNO]%
\renewcommand{\thesub}%
{(\arabic{section}.\arabic{THNO}.\arabic{sub})}%

\newcommand{\theab}%
{(\arabic{section}.\arabic{THNO}.\arabic{sub})}%

\newcounter{abcd}[sub]%
\renewcommand{\theabcd}%
{(\arabic{section}.\arabic{THNO}.\arabic{sub}.\alph{abcd})}%

\def\th#1{\refstepcounter{THNO}
\par\medskip\noindent\begingroup \it
{\bf\theTHNO\ #1\ }}%
\def\eth{\par\endgroup}

\def\subth#1{\refstepcounter{sub}%
\par\smallskip\noindent\begingroup\it {\bf\thesub\ #1\ }}%
\def\esubth{\par\endgroup}%

\def\de#1{\refstepcounter{THNO}%
\par\medskip\noindent\begingroup%
{\sc\theTHNO\ #1\ }}%
\def\ede{\par\endgroup}%

\def\subde#1{\refstepcounter{sub}%
\par\smallskip\noindent\begingroup{\sc\thesub\ #1\ }}%
\def\esubde{\par\endgroup}%

\def\ab{\refstepcounter{sub}%
\begingroup{\theab}\endgroup}%

\def\abc{\refstepcounter{abcd}%
\begingroup{\theabcd}\endgroup}%

\author{Yuri G. Prokhorov\thanks{Partially supported by the
Russian Foundation of Fundamental Research (??) and
the International Science Foundation (M 9o300).}}
\title{On the general elephant conjecture for Mori conic bundles}
\date{}
\begin{document}
\maketitle

\abstract{
Let $f:X\to S$ be an extremal contraction from a threefolds with terminal 
singularities onto a surface (so called Mori conic bundle). 
We study some particular cases of such 
contractions: quotients of usual conic bundles and index two contractions.
Assuming Reid's general elephants 
conjecture we also obtain a rough classification. We present many examples.
}

\section*{Introduction.}
This paper continues our study of extremal contractions from  threefolds
to  surfaces \cite{Pro}, \cite{Pro1}, \cite{Pro2}, \cite{Pro3}.
Such contractions occur naturally in birational classification theories
of three-dimensional algebraic varieties.
We are interested in the local situation.
\de{Definition.}
Let $(X,C)$ is a germ of three-dimensional normal complex space $X$ along
a compact reduced curve $C$ and let $(S,o)$ be a germ of a 
normal
two-dimensional complex space $S$ in  $o\in S$. Assume that 
$X$ has at 
worst
terminal singularities. We say that proper morphism 
$f:(X,C)\to (S,o)$ is a
{\it Mori conic bundle} if
\begin{itemize}
\item[{\rm (i)}]
$f^{-1}(o)_{\mt{red}}=C$,
\item[{\rm (ii)}]
$f_*\OOO_X=\OOO_S$,
\item[{\rm (iii)}]
$-K_X$ is $f$-ample.
\end{itemize}
The first example of Mori conic bundles are conic bundles in
the classical sense:\quad $f:(X,C)\to (S,o)$ is called
a (usual) {\it conic bundle} if $(S,o)$ is nonsingular and 
there exists an
embedding $i:(X,C)\hookrightarrow \PP^2\times(S,o)$ such that
$\OOO_{\PP^2\times S}(X)=\OOO_{\PP^2\times S}(2,0)$.
\ede
The general question arising very naturally lies in the 
classification
problem of Mori conic bundles.
The following two conjectures are interesting for 
application of the Sarkisov
program to study of birational properties of threefolds 
with structure of
conic bundle (see \cite{Iskovskikh}, \cite{Iskovskikh1}).

\th{Conjecture  I (special case of Reid's general elephants 
conjecture).}
\label{elephant}
Let $f: (X,C)\to (S,o)$ be a Mori conic bundle. Then 
a general member
of the anticanonical linear system $|-K_X|$ has only 
DuVal singularities.
\eth

\th{Conjecture II.}
 Let $f:(X,C)\to (S,o)$ be a Mori conic bundle. 
Then $(S,o)$ is either nonsingular
or a DuVal singularity of type $A_n$.
\eth
One can expect that a Mori conic bundle with very 
singular base is
a quotient of a usual conic bundle by a cyclic group (for example it 
follows
from the conjecture \ref{elephant}, see \ref{g.e.}).
So this particular case seems to be general. Such quotients 
are classified
in Section \ref{qb}. Surpisingly situation here is not very 
complicated and
there are only three cases. We also study index two Mori 
conic bundles
in Section \ref{s-ind-2}. Assuming conjecture \ref{elephant} 
we obtain a 
rough classification
of Mori conic bundles in Section \ref{sect-el}. In particular, 
we will show
(see also \cite{Pro}) that the first
conjecture implies the second one.  We present many examples.
The methods are completely elementary.
None of the techniques of Mori \cite{Mori-flip} or Shokurov
\cite{Shokurov} will be used.
Almost all the results of this work were announced in 
\cite{Pro1} and
\cite{Pro2}.
\par
{\sc Acknowledgments.}
I have been working on this problem at Max-Plank Institut 
f\"ur Matematik
in February-April 1995 and at the University of Warwick
in November-December 1995. The author would like to thank  staffs
of these institutes for hospitality. Different aspects of this work were
discussed with  Professors V.~A.~Iskovskikh, V.~V.~Shokurov and M.~Reid.
I am very grateful to them for help and advices.

\section{Preliminary results}
Throughout this paper a variety is a reduced irreducible complex space.
On a normal variety $X$ by $K_X$ we will denote its canonical (Weil) 
divisor.
If $X$ is a variety and $C\subset X$ is its closed subvariety, then
$(X,C)$ is a germ of $X$ along $C$. Sometimes we will replace
$(X,C)$ by its sufficiently small representative $X$.
When we say that  a variety $X$ has terminal singularities, it means
that singularities are not worse then that, so $X$ can be nonsingular.
\par
We need some facts about three-dimensional terminal singularities.

\de{}\label{terminal}
Let $(X,P)$ be a terminal singularity of index $m\ge 1$
and let $\pi:(X\3,P\3 )\to (X,P)$ be the canonical cover.
Then $(X\3,P\3 )$ is a terminal singularity of index 1.
It is known \cite{Pagoda} that $(X\3,P\3)$ is a
hypersurface singularity, i.~e. there exists an
$\cyc{m}$-equivariant embedding
$(X\3,P\3 )\subset (\CC^4,0)$.
\ede

\subth{Theorem \cite{Danilov}.}
\label{Danilov}
If in notations above $(X\3,P\3 )$ is smooth, then
it is $\cyc{m}$-isomorphic to $(\CC^3_{x_1,x_2,x_3},0)$ with the action of
$\cyc{m}$ by
$$
(x_1,x_2,x_3)\longrightarrow
(\varepsilon^{a} x_1,\varepsilon^{-a} x_2,\varepsilon^{b} x_3),
$$
where $\varepsilon=\exp(2\pi i/m)$,
and $a$, $b$ are integers prime to $m$.
Conversely every such singularity is terminal.
\esubth
Such singularity is denoted by $\frac{1}{m}(a,-a,b)$ or
$\CC^3/\cyc{m}(a,-a,b)$.
\subth{Theorem \cite{Mori-term}, \cite{RYPG}.}
\label{cl-term}
In notations above assume that $(X\3,P\3)$ is singular, then
it is $\cyc{m}$-isomorphic to a hypersurface
 $\{\phi (x_1,x_2,x_3,x_4)=0\}$ in $(\CC^4_{x_1,x_2,x_3,x_4},0)$,
where there are two cases for the action of $\cyc{m}$
\par
{\rm (i)} (the main series)
$(x_1,x_2,x_3,x_4;\phi)\longrightarrow
(\varepsilon^{a} x_1,\varepsilon^{-a} x_2,\varepsilon^{b} x_3,x_4;\phi),
$
where $\varepsilon=\exp(2\pi i/m)$,
and $a$, $b$ are integer prime to $m$.
\par
{\rm (ii)} (the exceptional series)  $m=4$ and
$(x_1,x_2,x_3,x_4;\phi)\to (ix_1,-ix_2,i^{a} x_3,-x_4;-\phi)$,
where $a=1$ or $3$.
\esubth
\subde{Remark.}
Terminal singularities of index $>1$ are classified by type
of general member $\in |-K_X|$. For example there are four
types of index two terminal singularities
$cA/2$, $cAx/2$, $cD/2$ and $cE/2$ \cite{RYPG}
(we use notations of \cite{KoM}).
\esubde

\th{Lemma.}
\label{st}
Let $(X,P)$ be a germ of terminal threefold singularity and let 
$F\in |-K_{(X,P)}|$
 be an element of anticanonical linear system. If $F$ is an
irreducible nonsingular surface, then $(X,P)$ is nonsingular.
\eth
\subde{Proof.}
If $(X,P)$ is of index one, then $F$ is Cartier and $(X,P)$ is nonsingular
in this case. So we assume that $(X,P)$ is of index $m>1$.
Consider the canonical cover $\pi:(X\3,P\3)\to (X,P)$ and let
$F\3:=\pi^{-1}(F)$. Since $\pi$ is \'etale outside $P$ and
$F-\{P\}$ is simply connected, the restriction $\pi:F\3\to F$
splits non-trivially, so $F\3=F\3_1+\dots+F\3_m$.
But then $P\3$ is the only point of intersection
of components $F\3_i$. On the other hand each $F\3_i$ is $\QQ$-Cartier,
a contradiction.
\qq
\esubde
It is well known that every DuVal  (and, more general,
log-terminal) singularity $(F,P)$ is
a quotient of a nonsingular germ $(\CC^2,0)$ by a finite group $G$
acting on $\CC^2$ free outside $0$. The order of $G$ is called
{\it the topological index } of $(F,P)$. Similar to \ref{st} one can prove
the following.
\subth{Lemma.}\label{index}
Let $(X,P)$ be a germ of a terminal threefold
singularity of index $m>1$ and  $F\in |-K_{(X,P)}|$
 be an anticanonical divisor.
Assume that the surface $F$ is reduced, irreducible and
the point $(F,P)$ is DuVal of topological index $n$.
Then $n$ is divisible by $m$. Moreover
 if $n=m$, then $(X,P)$ is a cyclic
quotient singularity and $(F,P)$ is of type $A_{m-1}$.
\esubth
Now we present some elementary properties of  Mori conic bundles.

\th{Theorem \cite{Cut} (see also \ref{non-sing}).}
\label{Cut}
Let $f:X\to S$ be a Mori conic bundle.
Assume that $X$ has only singularities of index $1$. Then
$S$ is a nonsingular surface and $f$ is a conic bundle (possibly
singular).
\eth
The following statement is a consequence of the Kawamata-Viehweg
vanishing theorem.
\th{Proposition (cf. \cite[(1.2)]{Mori-flip}).}
Let $f:X\to S$ be a Mori conic bundle.
Then
$$
R^if_*\OOO_X=0, \quad i>0.
$$
\eth

\subth{Corollary (cf. \cite[(1.3)]{Mori-flip}).}
Let $f:(X,C)\to (S,o)$ be a Mori conic bundle.
Then
\par
{\rm (i)} the fiber $C$ is a tree of rational curves,
i.~e. $p_a(C_0)=0$ for any one-dimensional subscheme $C_0\subset C$.
\par
{\rm (ii)}  \label{Pic}
$\mt{Pic}(X)\simeq H^2(C,\ZZ)\simeq\ZZ^{\rho}$, \quad
where $\rho$ is the number of components of $C$.
\esubth

\de{Construction.}
\label{construction}
Let $f:(X,C)\to (S,o)$  be a Mori conic bundle.
Assume that $(S,o)$ is a singular point.
By easy remark in \cite{KoMM} (see also \cite{Ishii}),
$(S,o)$ is a quotient singularity. It means that $(S,o)$ is a
quotient of nonsingular point $(S',o')$ by a finite group $G$,
where the action of $G$ on $S'-o'$ is free.     Therefore
there exists a faithful representation
$G\hookrightarrow GL(T_{o',S'})=GL_2(\CC)$.
Let $g:S'\to S$ be the quotient morphism and let $X'$ be the
normalization of $X\times_S S'$. We have the following commutative diagram
$$
\label{diagram}
\begin{array}{ccc}
(X',C')&\stackrel{h}{\longrightarrow}&(X,C)\\
\downarrow\lefteqn{\scriptstyle{f'}}&&\downarrow\lefteqn{\scriptstyle{f}}\\
(S',o')&\stackrel{g}{\longrightarrow}&(S,o)\\
\end{array}
$$
where $C':=f'^{-1}(o')$.
Then $G$ acts on $(X',C')$ and obviously $X=X'/G$, $C=C'/G$.
Since the action of $G$ on $S'-o'$ is free, so is the action
on $X'-C'$. Therefore $X'$ has only terminal singularities,
because $\mt{codim}(C')=2$ (see e.~g. \cite[6.7]{CKM}).
It gives us also that the action is free outside a finite number of
points $Q_1,\dots,Q_k\in X'$ (for each $Q_i$ its image $P_i=h(Q_i)$
has index $>1$).  Since $h:X'\to X$ has no ramification divisors,
the anticanonical divisor $-K_{X'}=h^*(-K_X)$ is ample over $S'$.
We obtain a new Mori conic bundle $f':(X',C')\to (S',o')$ over
a nonsingular base.
\par
Professor S. Mukai pointed out that actually $X\times_S S'$ is normal, so
we can take $X'=X\times_S S'$.
\ede
\th{Lemma (see \cite{Pro}, \cite{Pro1}, \cite{Ko}).}
\label{cyc}
In notations above
$G$ is cyclic and has at least one fixed point on $X'$.
In particular, $(S,o)$ is a cyclic quotient singularity.
\eth
\de{Proof.}
Since $p_a(C')=0$, it is easy to prove by induction on the
number of components of $C$ that $G$ has
either a fixed point $Q\in C'$ or  an invariant component $C_0'\subset C'$,
$C_0'\simeq \PP^1$.  But in the second case we have two inclusions
$G\subset PGL_2(\CC)$
and  $G\subset GL_2(\CC)$, where $G\subset GL_2(\CC)$ contains no
quasireflections. By the classification of finite subgroups in
 $PGL_2(\CC)\simeq SO_3(\CC)$, $G$ is cyclic and has two fixed points on
$C_0'$. Therefore in any case $G$ has a fixed point $Q\in X'$.
Take a small neighborhood $U\subset X$ of $P:=h(Q)$. We have a
surjective map $\pi_1(U-P)\to G$. But $\pi_1(U-P)$ is cyclic because
$(X,P)$ is a quotient of a hypersurface singularity by a cyclic group.
Whence so is $G$. This proves our claim.
\ede
\subth{Corollary (cf. \cite{Cut}).}
\par
{\rm (i)}
\label{non-sing}
If $X$ has index one, then $(S,o)$ is nonsingular.
\par
{\rm (ii) }
\label{|G|=2}
If $X$ has index two, then $(S,o)$ is either nonsingular or
DuVal of type $A_1$.
\par
{\rm (iii) }
\label{>n}
Let $f:(X,C)\to (S,o)$ be a Mori conic bundle.
If $(S,o)$ is a cyclic quotient singularity of type
$\frac{1}{n}(a,b)$, then $X$ contains at least one point of index
$\ge n$.
\esubth

\section{Quotients of conic bundles. }
\label{qb}
By  \ref{construction} and by  \ref{cyc} every Mori conic
bundle
$f:(X,C)\to (S,o)$ over a singular base is a quotient of
 another Mori conic bundle $f':(X',C')\to (S',o')$ with a nonsingular base
by a cyclic group. In this section we classify
Mori conic bundles $f:(X,C)\to (S,o)$ under the assumption that $X'$ is
Gorenstein (and then  $f':(X',C')\to (S',o')$ is a conic bundle by 
\ref{Cut}).
First we present several  examples.

\de{Example}\label{ex1}  (Toric example).
Let  $\PP^1 \times \CC^2 \to\CC^2$
be the standard projection. Define the action
of the group  $\cyc{n}$  on
$\CC^2_{u,v}$     and
$\PP^1_{x_0,x_1} \times\CC^2_{u,v}$:
$$
(x_0,x_1;u,v)\to (x_0,\varepsilon x_1;
\varepsilon^{a} u, \varepsilon^{-a}v),
$$
where $\varepsilon =\exp (2\pi i/n)$, $a\in\NN$ and $(n,a)=1$. Denote
$X=(\PP^1 \times\CC^2 )/\cyc{n}$,
$S=\CC^2 /\cyc{n}$.
Then the projection
 $f:X\to S$ is a Mori conic bundle.
The threefold  $X$ has on the fiber $f^{-1}(0)$
exactly two terminal points $P_1$, $P_2$
which are cyclic quotients of type
$\frac{1}{n}(a,-a,\pm 1)$, the surface
$S$ has in  0 a DuVal point of type $A_{n-1}$.
\ede

\de{Example.}\label{ex3}
Consider the following  hypersurface  in
$\PP^2_{x_0,x_1,x_2}\times\CC^2_{u,v}$:
$$
X':=\{x_0^2+vx_1^2+\psi(u,v)x_1x_2+ux^2_2=0\},
$$
where $\psi(0,0)=0$.
Let $n$ be an odd integer, $n=2q+1$, where $q\in\NN$. Define an action of
$\cyc{n}$ on $\PP^2\times\CC^2$ by
$$
(x_0,x_1,x_2,u,v)\to
(x_0,\varepsilon^{-q}x_1,\varepsilon^{q} x_2,\varepsilon u,\varepsilon^{-1}=
 v),
$$
where  $\varepsilon=\exp(2\pi i/n)$. If $\psi(u,v)$ is an
invariant, then $\cyc{n}$ acts naturally on $X'$.
As in \ref{ex1} the $f:X'/\cyc{n}\to\CC^2/\cyc{n}$ is a Mori conic bundle.
The singular locus of $X'/\cyc{n}$ consist of
two terminal cyclic quotient points of index $n$.
The point $(S,o)$ is DuVal of type $A_{n-1}$.
\ede

\de{Example.}\label{ex2}  Let $X'$ be a hypersurface  in
$\PP^2 _{x_0,x_1,x_2} \times\CC^2_{u,v}$, defined by the equation
$$
x_0^2+x_1^2+x_2^2\phi(u,v)=0,
$$
where $\phi(u,v)$ has no multiple factors and contains only monomials of
even degree.
Denote by $f':X'\to\CC^2$ the natural projection.
Then $X'$ has only one singular point
$P'=(x_0=x_1=u=v=0)$ on ${f'}^{-1}(0)$.
Define the action of $\cyc{2}$ on $X'$ and $\CC^2$ by
$$
(x_0,x_1,x_2, u,v)\to (-x_0,x_1,x_2,-u,-v).
$$
Let $X=X'/\cyc{2}$, $S=\CC^2/\cyc{2}$. The only fixed point on 
 $X'$ is $P'$
it gives us  a unique point $P\in X$ of index two.
The variety $X$ has no other singular points.
The surface $S$ has a DuVal singularity of type $A_1$ at $0$.
There are two cases for $\phi (u,v)$:
\par
(1) $\mt{mult}_{(0,0)}(\phi)=2$, then $(X,P)$ is terminal of 
type $cA/2$;
\par
(2) $\mt{mult}_{(0,0)}(\phi)\ge 4$, then $(X,P)$ is terminal of 
type $cAx/2$.
\ede

\th{Theorem.}
\label{th1}
Let $f:(X,C)\to (S,o)$ be a Mori conic bundle.
Assume that $f$ is a quotient of a conic bundle $f':(X',C')\to (\CC^2,0)$
by $\cyc{n}$ such that the action of $\cyc{n}$ on $\CC^2-\{0\}$ is free.
Then there exists an
analytic isomorphism between $f:(X,C)\to (S,o)$ and one of examples
\ref{ex1}, \ref{ex3} or \ref{ex2}.  In particular, $(S,o)$ is DuVal of type
$A_{n-1}$.
\eth
\de{Proof.}
By Cartan's lemma \cite{Car}, one can choose coordinates $u, v$ in
$\CC^2$ such that the action of $\cyc{m}$ is
$(u,v)\longrightarrow (\varepsilon^a u,\varepsilon^b v)$, where
$\varepsilon:=\exp (2\pi /m)$, $(a,m)=(b,m)=1$.
\par
Let $X'_0:=f'^{-1}(0)$ be the  scheme-theoretical fiber over $0$.
Then $X'_0$ isomorphic to a conic in $\PP^2$ (see \ref{Cut}).
There are the following cases.
\ede
\de{Case I.}  \underline{$X'_0$ is a non-degenerate conic.}
\label{caI}
Then in the analytic situation  $X'\simeq\PP^1\times\CC^2$
by the Grauert-Fisher theorem and  we may assume that the action of
$\cyc{n}$
in some coordinate systems $(x_0,x_1)$ in $\PP^1$ and $(u,v)$ in $\CC^2$ is
$$ (x_0,x_1;u,v)\longrightarrow 
(x_0,\varepsilon x_1;\varepsilon^a u,\varepsilon^b v),
\qquad\qquad (a,m)=(b,m)=1,
$$
where we may take $\wt(x_0)=0$ because $(x_0,x_1)$ are homogeneous
coordinates.
There are exactly two fixed points on $X'$:
$$
Q_0=\{x_0=u=v=0\},\qquad Q_1=\{x_1=u=v=0\}.
$$
They give us two points of index $n$ on $X'/\cyc{n}$ of types
$\frac{1}{n}(-1,a,b)$ and $\frac{1}{n}(1,a,b)$, respectively.
By \ref{Danilov} these two points  are terminal only if $a+b=n$
(recall that $a$, $b$ are defined modulo $n$.
We obtain the  example \ref{ex1}.
\ede

\de{Case II. } \underline{The fiber $X'_0$  is reducible.}
Then $X'_0=L_1+L_2$ is pair of lines intersecting each other in
one point, say $Q$, which is a fixed point. Let $P:=h(Q)$.
 We fix a generator $s\in\cyc{n}$ and for
$\cyc{n}$-semi-invariant  $z$ define weight $\wt(z)$ as an integer
defined modulo $n$ such that
$$
\wt(z)\equiv a\mod m\qquad  {\rm iff}
\qquad
s(z)=\varepsilon^az,
$$
where $\varepsilon =\exp{2\pi i/n}$.
\par
Similar to \ref{caI}
consider an $\cyc{n}$-equivariant
embedding $X'\subset \PP^2_{x_0,x_1,x_2}\times\CC^2_{u,v}$
such that $x_0,x_1,x_2$ are  semi-invariants with
$$
\wt(x_0,x_1,x_2;u,v)=(0,p,q;a,b),\qquad (a,n)=(b,n)=1,
$$
where we may take $\wt(x_0)=0$ because $(x_0,x_1,x_2)$ are homogeneous
coordinates. There are two subcases.
\subde{Subcase (i).} \underline{Components of $X_0'$ are 
$\cyc{n}$-invariant.}
\label{non-per}
We will derive a contradiction.
One can change homogeneous coordinates $x_0,x_1,x_2$ in $\PP^2$ so
that $X'_0\subset\PP^2$ is $\{x_0x_1=0\}$.
Then $X'$ is given by the equation
$$
x_0x_1+\varphi_0(u,v)x_0^2+\varphi_1(u,v)x_0x_1+\varphi_2(u,v)x_0x_2+
\phi(u,v)x_1^2+\psi(u,v)x_1x_2+\zeta(u,v)x^2_2=0,
$$
where $\varphi_0(0,0)= \varphi_1(0,0)= \varphi_2(0,0)= \phi(0,0)= 
\psi(0,0)= \zeta(0,0)=0$.
By taking $x_1'=x_1+\varphi_0x_0+\varphi_1x_1+\varphi_2x_2$
we obtain a new equation for  $X'$:
$$
x_0x_1+\phi(u,v)x_1^2+\psi(u,v)x_1x_2+\zeta(u,v)x^2_2=0,
\leqno
\abc
\label{equation}
$$
where $\phi(u,v)$, $\psi(u,v)$, $\zeta(u,v)$ are semi-invariants
with suitable weights and $\phi(0,0)=\psi(0,0)=\zeta(0,0)=0$.
Then $Q=\{x_0=x_1=u=v=0\}$   and
$y_0:=x_0/x_2, y_1=:x_1/x_2, u, v$ are local coordinates in
$\PP^2\times\CC^2$ near $Q$.
 We have
$$
(X,P)=\{y_0y_1+\phi(u,v)y_1^2+\psi(u,v)y_1+\zeta(u,v)=0\}
/\cyc{n}(-q,p-q,a,b).
\leqno
\abc
\label{pop}
$$
The action of $\cyc{n}$ on $X'$ has two more
fixed points $Q_1:=\{x_0=x_2=u=v=0\}$ and 
$Q_2:=\{x_1=x_2=u=v=0\}$.
Then
$z_0=x_0/x_1, z_2=x_2/x_1, u, v$ are local coordinates in
$\PP^2\times\CC^2$ near $Q_1$
and similarly
$t_1=x_1/x_0, t_2=x_2/x_0, u, v$ are local coordinates in
$\PP^2\times\CC^2$ near $Q_2$.
Let $P_i=h(Q_i)$, $i=1,2$.  Similar to \ref{pop}
$$
\begin{array}{l}
(X,P_1)=\{z_0+\phi(u,v)+\psi(u,v)z_2+\zeta(u,v)z_2^2=0\}
/\cyc{n}(-p,q-p,a,b)\simeq\\
\CC^3_{z_2,u,v}/\cyc{n}(q-p,a,b),\\
\end{array}
\leqno
\abc
\label{pop1}
$$
$$
\begin{array}{l}
(X,P_2)=\{t_1+\phi(u,v)t_1^2+\psi(u,v)t_1t_2+\zeta(u,v)t_2^2=0\}
/\cyc{n}(p,q,a,b)\simeq\\
\CC^3_{t_2,u,v}/\cyc{n}(q,a,b).
\end{array}
\leqno
\abc
\label{pop2}
$$
Since the action has a zero-dimension fixed locus,
$$
(q,n)=(p-q,n)=1.
\leqno
\abc
\label{=1}
$$
By \ref{cl-term} $(X,P)$ a cyclic quotient.
It is possible only if $\zeta(u,v)$ contains either
$u$ or $v$ terms. Up to permutation of $u, v$ we may assume that
$\zeta=u+\dots$. From \ref{equation} we have $\wt(x_0x_1)=\wt(ux_2^2)$.
Thus $p=a+2q$.
Then
$$
(X,P)\simeq \CC^3_{y_0,y_1,v}
/\cyc{n}(-q,p-q,b)=\CC^3/\cyc{n}(-q,a+q,b).
\leqno
\abc
\label{popp}
$$
We claim that $a+b=0$. Indeed otherwise $n>2$ and from \ref{pop1}
we have $b=a+q$ because $p=a+2q$ and $(q,n)=1$.
Whence $(X,P)=\CC^3/\cyc{n}(-q,b,b)$. This point cannot be terminal if 
$n>2$.
\par
The contradiction shows that $a+b=0$.
Point $(X,P)$ is terminal in this case only if $q=b$ (see \ref{popp}).
But then $p-q=a+2q-q=0$, a contradiction with \ref{=1}.
\esubde
\subde{Subcase (ii).} \underline{$\cyc{n}$ permutes components of $X_0'$.}
Then $n$ is even, $n=2k$. If $k>1$, then the quotient 
$X'/\cyc{k}\to\CC^2/\cyc{k}$
is a Mori conic bundle as in \ref{non-per}. We have proved that this is
impossible. Therefore $k=1$, $n=2$.
Then it is easy to see that $\wt(u)=\wt(v)=1$.
We may assume also that the fiber $X'_0\subset\PP^2$ is 
$\{x_0^2+x_1^2=0\}$.
Since $\cyc{2}$ permutes components of $\{x_0^2+x_1^2=0\}$,
we may assume that $\cyc{2}$ acts on $x_0$, $x_1$ by
$x_0\to -x_0$, $x_1\to x_1$.
Then one can change the coordinate system $(x_0,x_1,x_2;u,v)$ 
such that in
$\PP^2\times\CC^2$ the variety $X'$ is given by the equation
$$
x_0^2+x_1^2+\phi(u,v)x^2_2=0,
\leqno
\abc
\label{equation1}
$$
where $\phi(u,v)$ is an invariant with $\phi(0,0)=0$.
It means that $\phi$ contains only monomials of even degree.
Then $Q=\{x_0=x_1=u=v=0\}$ is the only fixed point on $C'$ and
$y_0:=x_0/x_2, y_1=:x_1/x_2, u, v$ are local coordinates in
$\PP^2\times\CC^2$ near $Q$. Whence
$$
(X,P)=\{y_0^2+y_1^2+\phi(u,v)=0\}/\cyc{2}.
$$
If this point is terminal, then up to permutation of $y_0$, $y_1$ 
we may
assume that $\wt(y_0)=1$, $\wt(y_1)=0$. Thus $\wt(x_2)=0$. Finally 
singularities
of $X'$ are isolated only if $\phi=0$ is a reduced curve.
Thus we have the example \ref{ex2}.
\esubde
\ede
\de{Case III.} \underline{$X'_0$ is a double line.}
As above  consider an $\cyc{n}$-equivariant
embedding $X'\subset \PP^2_{x_0,x_1,x_2}\times\CC^2_{u,v}$
such that the action of $\cyc{n}$ is
$$
\mt{wt}(x_0,x_1,x_2;u,v)=(0,p,q;a,b),\qquad (a,n)=(b,n)=1.
$$
We may assume also that the fiber $X_0'\subset\PP^2$ is $\{x_0^2=0\}$.

Then in some semi-invariant coordinate system $(x_0,x_1,x_2;u,v)$ in
$\PP^2\times\CC^2$ the variety $X$ is given by the equation
$$
x_0^2+\phi(u,v)x_1^2+\psi(u,v)x_1x_2+\zeta(u,v)x^2_2=0,
\leqno
\ab
$$
where $\phi(u,v)$, $\psi(u,v)$, $\zeta(u,v)$ are semi-invariants
with $\phi(0,0)=\psi(0,0)=\zeta(0,0)=0$.
As above we can take $\wt(x_0)=0$. Denote $p=\wt(x_1)$, $q=\wt(x_2)$.
The local coordinate along the fiber $X_0'$ is $x_1/x_2$ and
$\wt(x_1/x_2)=p-q$.
Since the action of $\cyc{n}$ on $X'$ is free in codimension two,
$(p-q,n)=1$.
Changing a generator of $\cyc{n}$ we can get $p-q=1$.
Fixed points are
$Q_1=\{ u=v=0, x_0=x_1=0\}$ and $Q_2=\{u=v=0, 
x_0=x_2=0\}$. We can take
the local coordinates near $Q_1$ and $Q_2$ in $\PP^2\times\CC^2$ as
and $(y_0=x_0/x_2,y_1=x_1/x_2,u,v)$ and 
$(z_0=x_0/x_1,z_2=x_2/x_1,u,v)$,
respectively.
The point $Q_1\in X'$ gives the singular point $P_1\in X$ of type
$$
\{y_0^2+\phi(u,v)y_1^2+\psi(u,v)y_1+\zeta(u,v)=0\}/\cyc{n}(-q,1,a,b).
\leqno\abc
\label{t1}
$$
Similarly, the point $Q_2\in X'$ gives the singular 
point $P_2\in X$ of 
type
$$
\{z_0^2+\phi(u,v)+\psi(u,v)z_2+\zeta(u,v)z_2^2=0\}/\cyc{n}(-p,-1,a,b).
\leqno\abc
\label{t2}
$$
\subde{}
We claim that $(X,P_1)$ and $(X,P_2)$ are from the main series.
 Indeed assume for example that $(X,P_1)$ is from the special series.
Then $n=4$, $q=2$, $p=3$. For $P_2$ all the weights are prime 
to $n=4$.
By \ref{cl-term}, $P_2$ is a cyclic quotient singularity and $X'$ is 
nonsingular at
$Q_2$. It is possible only if $\phi(u,v)$ contains a linear term.
Then $\wt(\phi)=\wt(z_0^2)=2$. It gives us $a=\wt(u)=2$ or 
$b=\wt(v)=2$,
a contradiction.
\esubde
\subde{}
Now we claim that $X'$ is nonsingular at $Q_1$ and $Q_2$.
As above if $Q_1\in X'$ is singular, then by \ref{cl-term} we have
$q=0$, $p=1$.
Again by \ref{cl-term}, $P_2$ is a cyclic quotient singularity and 
$X'$ is 
nonsingular at
$Q_2$. It is possible only if   $\phi(u,v)$ contains a linear term.
Up to permutation $u, v$  we may assume that $\phi=u+\cdots$.
Since $\wt(\phi)=\wt(z_0^2)=-2$, it
gives us $a=\wt(u)=-2$. Therefore  $n$ is odd and
$$
(X,P_1)=\{y_0^2+\phi(u,v)y_1^2+\psi(u,v)y_1+\zeta(u,v)=0\}/
\cyc{n}(0,1,-2,b).
\leqno\abc
\label{1}
$$
$$
\begin{array}{l}
(X,P_2)=
\{z_0^2+(u+\cdots)+\psi(u,v)z_2+\zeta(u,v)z_2^2=0\}/\cyc{n}(-1,-1,-2,b)\\
\simeq \CC^3_{z_0,z_2,v}/\cyc{n}(-1,-1,b).
\end{array}
\leqno\abc
\label{2}
$$
Since $P_2$ is a terminal point and $n$ is odd, from \ref{2} we 
have $b=1$.
But on the other hand  from \ref{1} $1+b=0$ or $-2+b=0$, a 
contradiction.
\esubde
\subde{}
Therefore $X'$ is nonsingular at $Q_1$ and $Q_2$.
It follows from \ref{t1} and \ref{t2} that both $\phi(u,v)$ and 
$\zeta(u,v)$
contain linear terms. Up to permutation of $u,v$ we may assume that
$\zeta=u+\dots$. Whence $a=-2q$.
Moreover $(X,P_1)=\frac{1}{n}(-q,1,b)$.
If $\phi=u+\dots$, then from \ref{t2} $-2p=a=-2+a$, so $n=2$, 
$a=0$,
a contradiction.  Therefore $\phi=v+\dots$. It gives us $-2p=b=a-2$ 
and
$(X,P_2)=\frac{1}{n}(-p,-1,a)$. Since $b$ is even, $n$ is odd.
Further $a=-2q$, $p=q+1$, $b=-2q-2$. Thus
$$
(X,P_1)=\frac{1}{n}(-q,1,-2q-2),\qquad\qquad 
(X,P_2)=\frac{1}{n}(-q-1,-1,-2q)
\leqno
\abc
\label{ppppp}
$$
\esubde
\subde{} Now we show that $n=2q+1$. Indeed assume the opposite.
Then since $(X,P_2)$ is terminal, from  \ref{ppppp} we have
either $q+2=0$ or $3q+1=0$. But in the first case
$(X,P_1)=\frac{1}{n}(2,1,2)$. This point is terminal only if $n=3$,
$q=1$. In the second case   $(X,P_1)=\frac{1}{n}(-q,1,-2q-2)$
cannot be terminal.
Thus our claim is proved. In particular, we have
$a=-2q=1$, $b=-2q-2=-1$.
Finally in this case $\wt(x_0,x_1,x_2,u,v)=(0,-q,q;1,-1)$ and
by changing coordinates in $\CC^2$ by $u'=\zeta(u,v)$, $v'=\phi(u,v)$
we obtain
$$
X\simeq
\{x_0^2+vx_1^2+\psi(u,v)x_1x_2+ux^2_2=0\}
/\cyc{2q+1}(0,-q,q;1,-1).
$$
Therefore $X$ is as in \ref{ex3}.
This proves our theorem.
\qq
\esubde

\section{Index two Mori conic bundles.}
\label{s-ind-2}
In this section  index two Mori conic bundles will be investigated.
First we consider the case  when the base
of a Mori conic bundle $f:(X,C)\to (S,o)$ is nonsingular, i.~e.
$(S,o)\simeq(\CC^2,0)$.
We need some elementary results about extremal 
neighborhoods
\cite{Mori-flip},
\cite{KoM}. Note that if $f:(X,C)\to (S,o)$ is a Mori conic bundle
with reducible central fiber $C$, then the  Mori cone
$\overline{NE}(X/S)$ is generated by extremal rays, because
$-K_X$ is ample. Since we consider a germ $(X,C)$, every extremal ray
 is generated by a component of $C$. On the other hand the dimension of
 $\overline{NE}(X/S)$ is equal to the number of components of $C$
by \ref{Pic}. Therefore $\overline{NE}(X/S)$ is simplicial and generated
by classes of components of $C$. So every irreducible component
$C_i\subset C$ gives us (not necessary isolated) extremal
neighborhood in the sense of \cite{Mori-flip}.

\th{Proposition.}
\label{en}
Let $(X,C)$ be an extremal neighborhood (not necessary isolated).
Assume that $X$ has index two and let $P\in X$ be an index two point.
Then
\subth{}  \cite[(4.6)]{KoM}
$P$ is the only point of index two.
$C$ has at most three components they all pass through $P$ and
they do not intersect elsewhere.
\esubth
\subth{}\cite[(2.3.2)]{Mori-flip}
\label{fu}
For every component  $C_i\subset C$ one has $(-K_{X}\cdot C_i)=1/2$.
\esubth
\subth{} \cite[(7.3)]{Mori-flip}
A general member $(F,P)\in |-K_{(X,P)}|$ satisfies
$F\in |-K_{X}|$, $F\cap C=\{ P\}$ and $(-K_{(X,C_i)}\cdot C_i)=1/2$
for every component $C_i\subset C$.
\esubth
\eth

The first main result of this section  is the following.

\th{Theorem.}
\label{th3}
Let $f:(X,C)\to (S,o)$ be a Mori conic bundle over a nonsingular  base 
surface
$(S,o)\simeq (\CC^2,0)$. Assume that $X$ is of index two. Then
\esubth
\subth{}
\label{-h}
$X$ contains exactly one index two point $P$, the central fiber $C$ has 
at most four components
they all  pass through $P$ and they do not intersect elsewhere.
\esubth
\subth{}
\label{h}
There exists  a flat elliptic fibration $g:(Y,L)\to (S,o)$
where $(Y,L)$  is a germ of  threefold with only isolated Gorenstein  
terminal
singularities along a reducible curve  $L$ such that  
$L=(g^{-1}(o))_{\mt{red}}$ 
and  a general fiber  $Y_s:=g^{-1}(s)$, $s\in S$ is a nonsingular 
elliptic curve.
Further $K_Y$ is trivial along $L$ and $g$ factores as
$$
g:(Y,L)\stackrel{h}{\longrightarrow} 
(X,C)\stackrel{f}{\longrightarrow} (S,o),
$$
where $h$ is a quotient morphism by  $\cyc{2}$
\esubth
\subth{}
\label{h1}
The action of $\cyc{2}$ on $L$ does not interchange components.
The locus of $\cyc{2}$-fixed points on $Y$ consists of an isolated point
 $Q$ such that $h(Q)=P$ and a nonsingular divisor $D\not\ni Q$ such that
$D\cap \mt{Sing}(X)=\emptyset$.
All the components of $L$ are isomorphic to $\PP^1$.
They  all pass through $Q$ and they do not intersect elsewhere.
\esubth
\subth{}
\label{-}
In notations above we have the following cases for
the scheme-theoretical fiber $X_o=f^{-1}(o)$.
In each case $C_1,\dots,C_r\simeq\PP^1$ are irreducible components of $C$.
\par\noindent\abc\quad
\label{1-1-1-1}
  $X_0\equiv C_1+C_2+C_3+C_4$,
\par\noindent\abc\quad
\label{1-1-2}
  $X_0\equiv C_1+C_2+2C_3$,
\par\noindent\abc\quad
\label{1-3}
  $X_0\equiv C_1+3C_2$,
\par\noindent\abc\quad
\label{2-2}
  $X_0\equiv 2C_1+2C_2$,
\par\noindent\abc\quad
\label{.4}
  $X_0\equiv 4C_1$.
\esubth
\eth
\de{Remark.}
Conversely, let $g:(Y,L)\to (S,o)\simeq(\CC^2,0)$
be an elliptic fibration with an action of $\cyc{2}$ such
as in \ref{h}-\ref{h1}.
If the point $(Y,Q)/\cyc{2}$ is terminal,
then $f:(X,C):=(Y,L)/\cyc{2}\to (S,o)$ is a Mori conic bundle
of index two.
\ede
We prove our theorem  in several steps.
\th{Lemma.}
\label{flat}
Let $f:X\to S$ be a morphism from a normal threefold with only
terminal singularities onto a surface. Assume that all the fibers of $f$ 
are
connected and one-dimensional. Then $f:X\to S$ is flat.
\eth
\de{Proof.}
Terminal singularities are rational \cite[1-3-6]{KMM} and therefore
Cohen-Macaulay \cite{Ke}. Then $f$ is flat by \cite[23.1]{Mat}.\qq
\ede
\de{}
Thus $f:X\to S$ is flat.
Let $X_s:=f^{-1}(s)$ be a scheme-theoretical fiber over $s\in S$.
Then $(X_o)_{\mt{red}}=C$ and $X_o\equiv \sum n_iC_i$ for some 
$n_i\in\NN$.
Since $X_s\simeq\PP^1$ for general $s\in S$, $(-K_X\cdot X_s)=2$.
Thus we have
$$
2=(-K_X\cdot X_s)=(-K_X\cdot X_o)=\sum n_i(-K_X\cdot C_i).
\leqno
\ab
\label{for}
$$
It gives us $\sum n_i\le 4$ because $(-K_X\cdot C_i)\in\frac{1}{2}\ZZ$.
In particular, $C$ has at most four components.
If $C$ is reducible then by \ref{fu} $(-K_X\cdot C_i)=1/2$, so
$\sum n_i=4$ and for $X_0$ we have only possibilities as in \ref{1-1-1-1}
-- \ref{2-2}.
\ede
\th{Lemma.}
Let $P\in X$ be a point of index 2. Then every component $C_i\subset C$
contains $P$.
\eth
\de{Proof.}
Assume $C_j\not\ni P$ for some $j$.
Then In particular, $C$ is reducible.
 By \ref{for}  we have
$\sum_{C_i\ni P}n_i<\sum n_i\le 4$.  Let  $F\in |-K_{(X,P)}|$
be a general member. Proposition \ref{en} gives us $(F\cdot C_i)=1/2$ 
for
 all
components  $C_i$ passing  through  $P$.  On the other hand
$$
1/2\sum_{C_i\ni P} n_i=\sum_{C_i\ni P} n_i (F\cdot  C_i)=(F\cdot X_0)
=(F\cdot X_s)\in\NN
$$
Thus
$1/2\sum_{C_i\ni P} n_i=1$. Therefore  $(F\cdot X_0)=(F\cdot X_s)=1$,
i.~e. the map $f|_F:F\to S$ is bimeromorphic and finite. Since
$S$ is nonsingular $f|_F$ is an isomorphism. So $F$ is also nonsingular.
We derive the contradiction with \ref{st}.
\qq
\ede
\th{Lemma.}
\label{tor-free}
The Weil divisor class group  $\mt{Cl}(X)$ has no torsion.
\eth
\de{Proof.}
If $\xi\in\mt{Cl}(X)$ is a torsion, then $\xi$ defines a cyclic \'etale
in codimension 1 cover $X'\to X$ (see e.~g. \cite[(1.11)]{Mori-flip}).
The threefold $X'$ is normal and has terminal
singularities of indices $\le 2$ only.
It follows also that $X'\to X$ is \'etale in codimension 2.
Take the Stein factorization
$$
\begin{array}{ccc}
X'&\longrightarrow&X\\
\downarrow&&\downarrow\\
S'&\longrightarrow&S\\
\end{array}
$$
We obtain a new Mori conic bundle $X'\to S'$ and \'etale in codimension 
one
cover $S'\to S$.  But $S$ is nonsingular. The contradiction shows that
$\mt{Cl}(X)$ is torsion-free.
\qq
\ede

\th{Lemma.}
\label{unique}
$X$ contains a unique point $P$ of index two.
\eth
\de{Proof.}
If $C$ is reducible, then $P$ is a unique point of
index two on $X$ by \ref{en} and because $p_a(C)=0$.
So assume that
 $C\simeq\PP^1$ and let
$P_1,\dots,P_r\in (X,C)$ be all the points of index $2$.

\subde{Definition.}
Let $X$ be a normal variety and $\mt{Cl}(X)$ be its
Weil divisor class group. The subgroup of
$\mt{Cl}(X)$ consisting of Weil divisor classes
which are $\QQ$-Cartier is called by the semi-Cartier
divisor class group. We denote it by $\mt{Cl}^{sc}(X)$.
\esubde

\subth{Theorem \cite{Pagoda},\cite{Kawamata}.}
\label{clsc}
Let $(X,P)$ be a germ of 3-dimensional terminal singularity.
Then $\mt{Cl}^{sc}(X,P)\simeq\cyc{m}$ and it is
generated by  the class of $K_{(X,P)}$.
\esubth

We have the following natural exact sequence
 (see \cite[(1.8.1)]{Mori-flip})
$$
\begin{array}{ccccccccc}
0&\to&\mt{Pic}(X)&\to&\mt{Cl}^{sc}(X)&\to&
\bigoplus_{i=1}^{r}\mt{Cl}^{sc}(X,P_i)&\to&0,\\
 &   &\|         &   &               &   &\|
                  &   &  \\
 &   &\ZZ        &   &               &   &\bigoplus_{i=1}^{r}\cyc{2} 
       &   &  \\
\end{array}
\leqno
\ab
\label{unique1}
$$
where $\mt{Pic}(X)\simeq\ZZ$ by \ref{Pic} and
$\mt{Cl}^{sc}(X,P_i)\simeq\cyc{2}$ by  \ref{clsc}. 
From \ref{tor-free} we get
$\mt{Cl}^{sc}(X)\simeq\ZZ$. Hence $r=1$. Thus our lemma is proved.
\qq
\ede
Since $p_a(C)=0$, components of $C$ do not intersect outside $P$.
This proves \ref{-h}.
\subth{Corollary.}
If $C$ is irreducible, then $(-K_X\cdot C)=1/2$.
\esubth
\subde{Proof.}
Suppose that  $(-K_X\cdot C)>1/2$
From=20lemma \ref{tor-free}  we have $\mt{Cl}^{sc}(X)=\ZZ$.
Moreover by \ref{unique1}
$\mt{Pic}(X)\subset \mt{Cl}^{sc}(X)$ is a subgroup of
index 2. Let $R$ be the ample generator of $\mt{Cl}^{sc}(X)$.
Then  $(R\cdot C)=1/2$ and $-K_X=kR$, $k\in\ZZ$.
On the other hand from
\ref{for} we have  $(-K_X\cdot C)=1$, $X_o\equiv 2C$
or $(-K_X\cdot C)=2$, $X_o\equiv C$. Whence $-K_X=2R$ or $-K_X=4R$.
But then $-K_X$ is Cartier and $X$ has index one, a contradiction.
\qq
\esubde
Therefore $X_0\equiv 4C$, if $C$ is irreducible.
This proves \ref{-}
\de{Construction.}
\label{construction2}
Let $B_i$ be a disc that intersects $C_i$ transversally in a general point
and let $B=\sum B_i$.
 Since $(-K_X\cdot C_i)=1/2$ and $\mt{Pic}(X)=\ZZ^r$, $B\in |-2K_X|$.
Take a double cover $h:Y\to X$ with ramification divisor $B$.
Set $L:=(h^{-1}(C))_{\mt{red}}$ and  $D:= (h^{-1}(B))_{\mt{red}}$.
We have
$$
g:(Y,L)\stackrel{h}{\longrightarrow} (X,C)\stackrel{f}{\longrightarrow}
(S,o)\simeq (\CC^2,0).
$$
By our construction $Y$ is normal, hence
 $g:Y\to\CC^2$ has only connected fibers.
 In  a neighborhood of each singular points on $X$
 $h$ is \'etale in codimension 1, therefore $Y$ has only terminal singular
points (see e.~g. \cite[(3.1)]{Pagoda} or \cite[(6.7)]{CKM}).
We have the following equalities for Weil divisors on $Y$
$$
K_{Y}=h^*(K_X)+D,\qquad\qquad D=h^*(-K_X).
$$
It gives us $K_{Y}=0$ and $Y$ has only terminal Gorenstein singularities.
The morphism $g$ is flat by lemma \ref{flat}.
Therefore $g$ is an elliptic fibration.
We have proved \ref{h}. It is clear that the locus of $\cyc{2}$-fixed
points on $Y$ consists of $D$ and $h^{-1}(P)$. This proves \ref{h1}.
Our theorem is proved.
\qq
\ede
Using \cite[(6.2)]{Mori-flip}, \cite[(4.7)]{KoM} and \cite{Pro1}
one can improve results of \ref{-}.
\subth{Corollary.}
In notations and conditions of \ref{th3} we have
\par\noindent\abc\quad
If  $X_0\equiv C_1+C_2+C_3+C_4$, then $P$ is the only singular 
point of $X$,
and $(X,P)$ is of type $cA/2$.
\par\noindent\abc\quad
If  $X_0\equiv C_1+C_2+2C_3$, then $X$ may have one more singular 
point of index one
on $C_3$,
and $(X,P)$ is of type $cA/2$.
\par\noindent\abc\quad
If $X_0\equiv C_1+3C_2$, then $X$ may have one more 
singular point of index one on $C_3$.
\par\noindent\abc\quad
If $X_0\equiv 2C_1+2C_2$, then $X$ may have one more singular point of 
index one
on  each component $C_i$, $i=1, 2$.
\par\noindent\abc\quad
If $X_0\equiv 4C_1$, then $X$ may have  at most two more 
singular points of one
 index 1.
\esubth

Now we use the construction \ref{construction2} to get examples index two
Mori conic bundles.
\de{Examples.}
In the following examples a Mori conic bundle $f:(X,C)\to (\CC^2,0)$
is constructed as a quotient $X=Y/\cyc{2}\to\CC^2$, where
$Y\subset\PP^3\times\CC^2$ is an intersection of two quadrics,
$\cyc{2}$ acts on $Y\subset\PP^3_{x_0,x_1,x_2,x_3}\times\CC^2_{u,v}$ by
$$
(x_0,x_1,x_2,x_3;u,v)\longrightarrow 
(-x_0,-x_1,-x_2,x_3;u,v)
$$
and $f$ is induced by the projection of 
$g:Y\subset\PP^3\times\CC^2\to\CC^2$ on
the second factor. We will use the notations $cA/2$, $cAx/2$, $cD/2$,
$cE/2$ to distinguish index two terminal singularities \cite{KoM}.
\ede

\de{Example.}
Let $Y\subset\PP^3_{x_0,x_1,x_2,x_3}\times\CC^2_{u,v}$
is given by the equations
$$
\left\{ \begin{array}{l}
x_0x_1=(au+bu^2+cuv)x_3^2\\
(x_0+x_1+x_2)x_2=vx_3^2,
\end{array}\right.
$$
where $a, b, c \in \CC$ are constants.
It is easy to check that $Y$ is nonsingular near the central fiber
$g^{-1}(0)$ if $a\ne 0$ and has an isolated hypersurface singularity in
$u=v=x_0=x_1=x_2=0$ if $a=0$, $c\ne 0$, $b\ne 0$ such that the 
rank of the
quadratic part is equal to $3$.
Then the quotient $X:=Y/\cyc{2}\to \CC^2$ is a Mori conic bundle with
only one singular point of type $cA/2$ and the central fiber
$X_o\equiv C_1+C_2+C_3+C_4$.
If $a\ne 0$, then the singular point is a cyclic quotient of type
$\frac{1}{2}(1,1,1)$.
\ede

\de{Example.}
Let $Y\subset\PP^3_{x_0,x_1,x_2,x_3}\times\CC^2_{u,v}$
is given by the equations
$$
\left\{ \begin{array}{l}
x_0x_1=u(x_0^2+x_1^2+x_2^2-x_3^2)+avx_3^2\\
(x_0+x_1)x_2=vx_3^2+bux_0^2
\end{array}\right.
\qquad\qquad a, b \in \CC
$$
Then the quotient $X:=Y/\cyc{2}\to \CC^2$ is a Mori conic bundle with
the central fiber $X_o\equiv C_1+C_2+2C_3$.
$X$ contains exactly one non-Gorenstein point that is
of type $\frac{1}{2}(1,1,1)$. If $a=b=0$, then
$X$ has also an ordinary double point on $C_3$.
\ede

\de{Example.}
Let $Y\subset\PP^3_{x_0,x_1,x_2,x_3}\times\CC^2_{u,v}$
is given by the equations
$$
\left\{ \begin{array}{l}
x_0x_1-x_2^2=ux_3^2\\
x_0x_2=ux_1^2+v(x_2^2+x_3^2)
\end{array}\right.
$$
Then the quotient $X:=Y/\cyc{2}\to \CC^2$ is a Mori conic bundle with
only one singular point of type $\frac{1}{2}(1,1,1)$ and the central fiber
$X_o\equiv C_1+3C_2$.
\ede
The following example shows that all the types of terminal singularities
of index two can appear on Mori conic bundles such as in \ref{2-2}.

\de{Example.}
Let $Y\subset\PP^3_{x_0,x_1,x_2,x_3}\times\CC^2_{u,v}$
is given by the equations
$$
\left\{ \begin{array}{l}
x_0x_1=ux_3^2\\
x_2^2=u(x_0^2+x_1^2)+(av+bv^2+cv^3+duv+eu^2v)x_3^2\\
\end{array}\right.
$$
where $a,b,c,d,e\in\CC$ are constants. If at least one of
$a,b,c$ is non-zero, then $Y$ has near the central fiber $Y_0$ an
isolated singularity at $\{x_0=x_1=x_2=u=v=0\}$ (or nonsingular).
As above the quotient $X:=Y/\cyc{2}\to \CC^2$ is a Mori conic bundle with
the central fiber $X_0\equiv 2C_1+2C_2$ and only one singular point
$P:=C_1\cap C_2$.
We have the following possibilities:
\par
$a\ne 0$, then $Y$ is nonsingular, so $P\in X$ is of type
$\frac{1}{2}(1,1,1)$,
\par
$a=0$, $b\ne 0$, then $P\in X$ is of type $cAx/2$,
\par
$a=b=e=0$, $c\ne 0$, $d\ne 0$, then  $P\in X$ is of type $cD/2$,
\par
$a=b=d=0$, $c\ne 0$, $e\ne 0$, then  $P\in X$ is of type $cE/2$.
\ede

\de{Example.}
Let $Y\subset\PP^3_{x_0,x_1,x_2,x_3}\times\CC^2_{u,v}$
is given by the equations
$$
\left\{ \begin{array}{l}
x_0^2=ux_2^2+vx_3^2\\
x_1^2=ux_3^2+vx_2^2
\end{array}\right.
$$
Then the quotient $X:=Y/\cyc{2}\to \CC^2$ is a Mori conic bundle with
irreducible central fiber that  contains
one singular point of type $\frac{1}{2}(1,1,1)$ and two ordinary double 
points.
\ede
Now we investigate index two Mori conic bundles with singular base.

\th{Theorem.}
Let $f:(X,C)\to (S,o)$ be a Mori conic bundle of index 2 over a surface.
Assume that the point $(S,o)$ is singular. Then $(S,o)$  is DuVal
of type $A_1$ and $f$ is either as in example \ref{ex2} or as in example
\ref{ex1} with $n=2$.
\eth
\de{Proof.} We can use construction \ref{construction}.
Since  $X$ contains only points of indices $\le 2$ by \ref{>n},
$(S,o)$ is a singularity of type $\frac{1}{2}(1,1)=A_1$. We claim that
$f':(X',C')\to (\CC^2,0)$ is a conic bundle.
It is sufficient to show  that $X'$ is Gorenstein \ref{Cut}.
Assume the opposite.
We remark that $X'$ contains only
points of indices $\le 2$ because $X'\to X$ is \'etale in codimension 1.
 By theorem \ref{th3}
$X'$ contains a unique point, denote it by $Q$, of index 2.
Then $Q$ is a $\cyc{2}$-fixed point. But then the point $P:=h(Q)$
has index $4$, a contradiction. Therefore $f'$  is a conic bundle and
by theorem \ref{th1} $f(X,C)\to (X,C)$ is such as in  example \ref{ex1} or
example \ref{ex2}. This proves our theorem.
\qq
\ede

\section{The general elephant conjecture.}
\label{sect-el}
In this section we study  Mori conic bundles
 under the assumption the existence of a good member in $|-K_X|$.
\de{Example.}
\label{exc3}
Let $f:(X,C)\to (S,o)$ be as in example \ref{ex3}. Consider the open subset
$U_2=\{ x_2\ne 0\}\subset X$. The local coordinates in $U_2$
$(t_0=x_0/x_2, t_1=x_1/x_2, v)$. Consider also the rational 
differential
$\sigma=(1/t_0)(dt_0\wedge dt_1\wedge dv)$ on $U_2$. Then it is easy to 
see that
$\sigma$ can be extended on $X'$ near $C'$.
Since $\sigma$ is $\cyc{n}$-invariant, $\sigma^{-1}$ defines
an element $F\in |-K_X|$, the image of $\{x_0=0\}$. It is easy to check
that $F$ contains the central
fiber $C=(f^{-1}(o))_{\mt{red}}$ and has two singular points of 
type $A_{n-1}$.
\par
Similarly one can check that in example \ref{ex1} a general member
$F\in |-K_X|$ does not contain $C$ and has two connected components.
It contains two singular points of type $A_{n-1}$ on each component.
\ede
The following theorem improves results of \cite{Pro}.
\th{Theorem.}
\label{g.e.}
Let $f:(X,C)\to (S,o)$ be a Mori conic bundle. Assume that conjecture
\ref{elephant} holds. Then we have one of the following:
\subth{}
\label{part1}
$(S,o)$ is nonsingular,
\esubth
\subth{}
\label{part2}
$(S,o)$ is DuVal of type $A_1$,
\esubth
\subth{}
\label{imp}
$(S,o)$ is DuVal of type $A_3$, in this case $C$ is irreducible,
$(X,C)$ has a cyclic quotient
singularity $P$ of index $8$ and has no another points of index $>1$
(see \cite{Pro3} for more detailed study of this case).
\esubth
\subth{}
\label{pr}
$f:(X,C)\to (S,o)$ is quotient of a nonsingular conic bundle
$f':(X',C')\to (S',o')$ with irreducible $C'$ by the group $\cyc{n}$, where
$n\ge 3$ and the action $\cyc{n}$ on $(S',o')\simeq (\CC^2,0)$
is free in codimension 1 (i.~e. $f:(X,C)\to (S,o)$ is such as in
\ref{ex1} or \ref{ex3}). In particular, $(S,o)$ has type $A_{n-1}$
in this case.
\esubth
\eth
\de{Proof.}
\label{as}
Assume that $(S,o)$ is singular and it is not of type $A_1$.
We will use the notations of construction \ref{construction}.
If $X'$ is of index one, then by \ref{Cut} $(S,o)$ is nonsingular.
So we assume that $X'$ is of index $>1$ and show in this case that
$f:(X,C)\to (S,o)$ is such as in \ref{imp}.
 Let $F\in |-K_X|$ be a general member and $F':=h^{-1}(F)$.
By our conditions $F$ has only DuVal singularities and since
$F'\to F$ is \'etale in codimension 1, so has $F'$ (In particular, $F'$
is irreducible).
Since $(-K_X\cdot X_s)=2$, where $X_s$ is a general fiber of $f$, the 
restriction
$f|_F:F\to S$ is  generically finite of degree 2. There are the following
cases.
\ede
\de{Case I.} \underline{$F'\cap C'$ is disconnected}.
As above since $F'\in |-K_{X'}|$, we see that $f|_{F'}:F'\to S'$ is
generically finite of degree 2 and $F'$ has two connected components
$F'_1$ and $F'_2$. Let
$$
f':F'\stackrel{f'_1}\longrightarrow D'\stackrel{f'_2}\longrightarrow S'
$$
be the Stein factorization.
Then $f'_1:F'\to D'$ is bimeromorphic  and $f'_2:D'\to S'$ is finite of
degree 2. In our case $D'$ has exactly two irreducible components
$D'_1$ and $D'_2$. Therefore we have $D_1'\simeq D_2'\simeq S'$ and
the divisor $D'$ is nonsingular, because so is $S'$. On the other
hand by the adjunction formula, $K_{F'}=0$.
Whence the morphism $f_1'$ is crepant (i.~e. $K_{F'}={f_1}'^*K_{D'}$).
It means that $F'$ is nonsingular and $f'_1=\mt{id}$. By
\ref{st}, $X'$ has no points of indices $>1$, a contradiction with
our assumption in \ref{as}.

\ede
\de{Case II.} \underline{$F'\cap C'$ is one point $Q$}.
\label{finite}
Then $F\cap C=\{ h(Q)\}$ is also one point, say $P$.
Therefore the morphism $f|_F:F\to S$ is finite of degree 2 and
 $P$ is a unique point of index  $>1$ and $(F,P)\to (S,o)$ is
double cover of isolated singularities. Thus $(S,o)$ is a
quotient of DuVal singularity $(F,P)$ by an involution $\tau$.
\par
Actions of involutions on DuVal singularities were classified by Catanese
\cite{Cat}. Recall that $(S,o)$ is of topological index $n\ge 3$
(otherwise we have cases \ref{part1} or
\ref{part2} of our theorem). Taking into account that $(S,o)$ is a 
cyclic quotient from the
list in \cite{Cat} we obtain the following
\th{Lemma \cite{Cat}.}
\label{Cat}
Let $(F,P)$ be a germ of DuVal singularity and let $\tau$ be an analytic
involution acting on $(F,P)$. Assume that the quotient 
$(S,o):=(F,P)/\tau$
is a cyclic quotient singularity of type $\frac{1}{n}(a,b)$ with $n\ge 3$
(and $(a,n)=(b,n)=1$). Then there are the following possibilities
for $(F,P)\to (S,o)$:
$$
\begin{array}{lll}
&(F,P)\to (S,o)&n\\
&&\\
\ab\label{(1)}\qquad&E_6\stackrel{2:1}{\longrightarrow}A_2, &n=3,\\
\ab\label{(2)}\qquad&A_{2n-1}\stackrel{2:1}{\longrightarrow}A_{n-1},
&n\ge 1\\
\ab\label{(3)}\qquad&A_{2k}\stackrel{2:1}{\longrightarrow}
\frac{1}{2k+1}(k,2k-1),&n=2k+1\\
\ab\label{(4)}\qquad&A_k\stackrel{2:1}{\longrightarrow}A_{2k+1}, 
&n=2k+1,\\
\ab\label{(5)}\qquad&A_{2k+1}\stackrel{2:1}{\longrightarrow}
\frac{1}{4k+4}(2k+1,2k+1),&n=4k+4.\\
\end{array}
$$
\eth

The restriction $(F',Q)\to (F,P)$ is \'etale outside $P$ and of degree $n$.
Therefore $(F',Q)$ is DuVal. It is easy to see (see e.~g. 
\cite[4.10]{RYPG}),
that  $(F,P)$ has no such covers in cases \ref{(4)}, \ref{(5)}. The 
singularity
$(F,P)$ from \ref{(3)} admits only cover by nonsingular $(F',P')$ of degree
$n=2k+1$. But then $(X',Q)$ is a nonsingular point by \ref{st}, a
contradiction with our assumption in \ref{as}.
\par
Let $m$ be the index of $(X,P)$, $\pi:(X\3,P\3)\to (X,P)$ be the canonical
cover and $F\3:={\pi\3}^{-1}{F}$. As above  we have \'etale in 
codimension one
$\cyc{m}$-cover $\pi:(F\3,P\3)\to (F,P)$ of DuVal singularities.
Since $\pi$ factores through $(X',Q)$ we have $m\ge n$.
\par
In case \ref{(1)} $(F,P)=E_6$ by \cite[4.10]{RYPG}, admits only 
cyclic cover
$D_4\stackrel{3:1}{\longrightarrow}E_6$. Then $n=m=3$. Therefore
$(X\3,P\3)\simeq (X',P')$ has index 1, a contradiction.
\par
Finally, consider case \ref{(2)}. Since $(X',P')$ has index $>1$, we have
$m>n\ge 3$ and by \cite[4.10]{RYPG}
$(F\3,P\3)\stackrel{m:1}{\longrightarrow}(F,P)$ is      of type
$(\mt{nonsingular})\stackrel{2n:1}{\longrightarrow}A_{2k+1}$, $m=2n$.
Then the index of $(X',P')$ is equal to $m/n=2$, $(X\3,P\3)$ is 
nonsingular,
hence  $(X',P')$ is a cyclic quotient singularity of type
$\frac{1}{2}(1,1,1)$.
In particular, $f$ is as in theorem \ref{th3}.
 Let $X_{o'}'=f'^{-1}(o')$ be the scheme-theoretical fiber
of $f'$ over $o'$.  Further $P'$ is the only fixed point on $C'$ under
the action of $\cyc{n}$. Whence $\cyc{n}$ permutes components
$\{C_i'\}$ of $C'$. In particular, the number of components $\ge n\ge 3$
and multiplicities of components in $X_{o'}'$ are the same.
Therefore for $f':(X',C')\to (S',o')$ we have the only possibility
\ref{1-1-1-1}. Hence we have exactly four components
of $X_0'$, $n=4$ and the fiber $X_o'$ is reduced. Then $X'$ is 
nonsingular outside $P'$. We obtain case
\ref{imp} of our theorem.
\ede

\de{Case III.} \underline{$F'\cap C'$ is one-dimensional and connected}.
Then so is $F\cap C$. Let $L:=F'\cap C'$ (with reduced structure).
On the surface $F'=h^{-1}(F)$ we have by the adjunction formula
that $L$ is Gorenstein and $\omega_L=\OOO_L(L)$, because
$\omega_{F'}=\OOO_{F'}$.
Since $L$ is contracted by $f'$, the dualizing sheaf $\omega_L$ is 
anti-ample.
Whence $L$ is a reduced conic in $\PP^2$ (see e.~g. \cite{Lip}).
\par
Let
$$
f_F:F\stackrel{f_1}\longrightarrow D\stackrel{f_2}\longrightarrow S
$$
be the Stein factorization. Then $f_1:F\to D$ is bimeromorphic and
$f_2:D\to S$ is finite of degree 2. By the adjunction formula, $K_F=0$.
Therefore the morphism $f_1$ is crepant (i.~e. $K_F=f_1^*K_D$) and $D$ 
has only  DuVal singularities.
\par
There exists the common minimal resolution 
$\sigma:\widetilde{F}\to F\to D$.
In our case $f_2^{-1}(o)$ consist of only one point $R$. Let
$\Gamma=\Gamma(\widetilde{F}/D)$ be a dual graph for $\sigma$. Denote 
vertices
corresponding to components of $h(L)$ (resp. $f_2$-exceptional divisors) 
by $\bullet$ (resp.
$\circ$). Then white vertices form connected subgraphs corresponding 
singular
points of $(F,C)$ and black vertices correspond components of $C$ that
contained in $F$. For $f_2:(D,R)\to (S,o)$ we have the same possibilities
as for $(F,P)\to (S,o)$ in \ref{Cat}. Let us consider these cases.

\subde{Subcase \ref{(1)}.} (i.~e. $(D,R)=E_6$, $(S,o)=A_2$, $n=3$).
The group $\cyc{3}$ naturally acts on the curve $L=F'\cap C'$.
Since $L$ is a conic, $\cyc{3}$ cannot permute its components.
Therefore $\cyc{3}$ has two or three fixed points  $Q_i$ on $L$.
Then there exists at least two points  $P_1,P_2\in X$ of indices $\ge 3$.
These points on $F$ have  by \ref{index} topological indices
$\ge 3$. On the other hand, the graph of the minimal
resolution of $(F,P_i)$ must be a white subgraph of 
$\Gamma$ ($\simeq E_6$).
Whence each $(F,P_i)$ is of type $A_2$ and there are only two fixed points
on $L$. It is possible only $L$ is irreducible and so is $h(L)$.
Thus we have only one possibility for $\Gamma$.
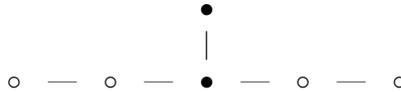
$$
\begin{array}{ccccccccc}
&&&&\bullet&&&&\\
&&&&|&&&&\\
\circ&\pal&\circ&\pal&\bullet&\pal&\circ&\pal&\circ\\
\end{array}
$$
Whence $F$ contains exactly two singular points $P_1$, $P_2$
and they have type $A_2$. The variety $X$ has at these points
index 3. Since $h^{-1}(P_i)=Q_i$ by \ref{index}, $X'$ is Gorenstein.
We obtain a contradiction with assumptions in \ref{as}.
\esubde

\subde{Subcase \ref{(2)}} i.~e.  $(D,R)=A_{2n-1}$, $(S,o)=A_{n-1}$.
If $L$ is irreducible or $\cyc{n}$ doesn't permute components of $L$,
then as above $F$ contains at least two points of topological indices 
$\ge n$.
Then $\Gamma$ is
$$
\underbrace{\circ\pal\circ\pal\cdots\pal\circ}_{n-1}
\pal\bullet\pal\underbrace{\circ\pal\cdots\pal\circ}_{n-1}
\leqno
\abc
\label{lll}
$$
Whence $F$ contains exactly two singular points and they have type
$A_{n-1}$. By \ref{index} $X'$ is Gorenstein, a contradiction.
\par
Therefore $L$ is reducible and $\cyc{n}$ interchanges its components.
Then $h(L)=F\cap C$ is irreducible, so $\Gamma$ has only
one black vertex. Further a $\cyc{n}$-fixed point $Q_1$ on $X'$
gives us the point $P_1\in F$ of topological index $nk$, $k\in\NN$.
This point corresponds to a white subgraph in $\Gamma$ of type
$A_{nk-1}$. Thus $\Gamma$ again has the form \ref{lll}
So $(F,P_1)$ is of type $A_{n-1}$.
In particular, by \ref{st} $Q_1$ is nonsingular.
The second white subgraph in $\Gamma$ also corresponds to the singular 
point
$P_2\in F$ of type $A_{n-1}$. Then the index of $(X,P_2)$ divides $n$.
Since $\cyc{n}$ interchanges two components of $L$, $n=2k$ is even
and there are three $\cyc{k}$-fixed points on $L$:
$Q_1$, $Q_2$ and $Q_3$, where $h(Q_2)=h(Q_3)=P_2$.
Whence the index of $P_2$ is divisible by $k$.
By our assumption $(X',C')$ is not Gorenstein, so
$(X',Q_2)$, $(X',Q_3)$ have
(the same) index $m_0>1$ and $(X',C')$
contains no another points of index $>1$.
But then index of $(X,R)$ is
$m_0n/2\le n$.  Hence $m_0=2$, it means that $f':(X',C')\to (S',o')$
is an index two Mori conic bundle that contains two non-Gorenstein points
$Q_2$, $Q_3$, a contradiction with \ref{th3}.
\esubde
Finally subcases \ref{(3)},\ref{(4)} and \ref{(5)} are impossible as above.
\qq
\ede

\section{An example of Mori conic bundle of index three.}
In this section we construct an example of index three Mori conic bundle
over a nonsingular base (that is not such as in examples \ref{ex1} or
\ref{ex3}).
\de{} Let $V\subset\PP^3_{x_0,x_1,x_2,x_3}\times\CC^2_{u,v}$ is given
by two equations
$$
\left\{
\begin{array}{l}
a_{1}x_{1}^{3}+a_{2}x_{2}^{3}+
a_{3}x_{1}^{2}x_{2}+a_{4}x_{2}^{2}x_{1}+
a_{5}x_{3}^{2}x_{1}+a_{6}x_{3}^{2}x_{2}=ux_{0}^{3}\\
b_{1}x_{1}^{3}+b_{2}x_{2}^{3}+
b_{3}x_{1}^{2}x_{2}+b_{4}x_{2}^{2}x_{1}+
b_{5}x_{3}^{2}x_{1}+b_{6}x_{3}^{2}x_{2}=vx_{0}^{3}\\
\end{array}
\right.
$$
where  $a_i,b_i\in\CC$ are general constants. Then $V$ is nonsingular
    and
the projection on the second factor $q:V\to\CC^2$ gives us a fibration
of curves of degree 9 in  $\PP^3_{x_0,x_1,x_2,x_3}\times\CC^2_{u,v}$.
 The central fiber $q^{-1}(0)$ has exactly nine
irreducible components
$\Gamma_0,\dots,\Gamma_8$, where $\Gamma_0:=\{u=v=x_1=x_2=0\}$.
Define an action of $\cyc{6}$ on $V$ by
$$
(x_0,x_1,x_2,x_3;u,v)\longrightarrow (\varepsilon x_0,x_1,x_2,-x_3;u,v),
\qquad \varepsilon =\exp (2\pi i/6).
$$
First we consider the quotient $p: V\to Y:=V/\cyc{3}$ and the standard 
projection
$g:Y\to\CC^2$. The set of $\cyc{3}$-fixed points is
a divisor $R:=\{x_0=0\}\cap V$ and an isolated point
$O:=\{x_1=x_2=x_3=u=v=0\}$.
The local coordinates on $V$ near $O$ are $(x_1/x_0, x_2/x_0, x_3/x_0)$
with $\cyc{3}$-weights  $\wt(x_1/x_0)=\wt(x_2/x_0)=\wt(x_3/x_0)$.
Therefore the singular locus of $Y$ consists of one point of type
$\frac{1}{3}(1,1,1)$ (denote it by $Q$) that is canonical and Gorenstein.
Further by the Hurvitz formula
$$
K_V=p^*K_Y+2R, \qquad  p^*K_Y=K_V-2R=0.
$$
Whence the canonical divisor of $Y$ is trivial and
$g:Y\to \CC^2$ is an elliptic fibration.
Since $\cyc{3}$ does not permute  $\Gamma_0,\dots,\Gamma_8$,
the central fiber $g^{-1}(0)$ consists of nine components 
$L_0,L_1,\dots,L_8$
which are proper transforms of corresponding components
$\Gamma_0,\dots,\Gamma_8$  of $q^{-1}(0)$.
Now consider the quotient $h:Y\to X:=Y/\cyc{2}$ and the natural morphism
$f:X\to\CC^2$.
The component $L_0$ is $\cyc{2}$-invariant and $\cyc{2}$ permutes
$L_i$, $i=1,\dots,8$ non-trivially.
 Hence $f^{-1}(0)=g^{-1}(0)/\cyc{2}$ has exactly five components.
Denote them by $C_0,\dots,C_4$.
To prove that $f:X\to\CC^2$ is  a Mori conic bundle, we have to investigate
the singular locus of $X$.
The set of $\cyc{2}$-fixed points on $V$ is an irreducible divisor 
$D:=\{x_3=0\}$.
Therefore the set of $\cyc{2}$-fixed points on $Y$ is an irreducible
 Weil divisor $F:=p(D)$ such that $3F$ is Cartier.
Moreover $D$ intersects components $\Gamma_0,\dots,\Gamma_8$
transversally.
\ede
The following is an easy exercise.
\th{Lemma.}
Let $(Y,Q)$ be a cyclic quotient singularity of type
$\frac{1}{3}(1,1,1)$ with action of $\cyc{2}$.
Assume that the locus of fixed points of this action
is a Weil divisor $F$ such that $3F=0$.
Then  $(Y,Q)/\cyc{2}$ is terminal of type $\frac{1}{3}(1,1,-1)$.
\qq
\eth

\de{}
By the Hurvitz formula
$0=K_Y=h^*(K_X)+F$. Since $(L_i\cdot F)>0$,
we have $(h^*(-K_X)\cdot L_i)=(F\cdot L_i)>0$, $i=0,\dots,8$.
Hence $(-K_X\cdot C_i)>0$, $i=0,\dots,4$ i.~e.
$-K_X$ is relatively $f$-ample (at least near $f^{-1}(0)$.
\ede
\de{Conclusion.}
We get   a Mori conic bundle $f:(X,C)\to (\PP^2,0)$ such that
its central fiber $C$ has exactly five components
and the only singular point $P$ of $X$ is a cyclic quotient singularity of
type $\frac{1}{3}(1,1,-1)$. All the components of $C$ pass through
$P$ and they do not intersect elsewhere. Note that $h:Y\to X$ is nothing
else but Kawamata's double cover trick \cite{Kawamata}.
\ede

\par\noindent
{\sc Algebra Section,
Department of Mathematics,
Moscow State University,
Vorob'evy Gory, Moscow
117 234, Russia}

\par\noindent
\begin{tabular}{ll}
{\sc E-mail:}&prokhoro@mech.math.msu.su\\
&prokhoro@nw.math.msu.su\\
\end{tabular}

\end{document}